\documentclass[final,onefignum,onetabnum]{siamart171218}

\usepackage{lipsum}
\usepackage{amsfonts}
\usepackage{graphicx}
\usepackage{epstopdf}
\usepackage{algorithmic}
\usepackage{multirow}

\usepackage[english]{babel}
\ifpdf
  \DeclareGraphicsExtensions{.eps,.pdf,.png,.jpg}
\else
  \DeclareGraphicsExtensions{.eps}
\fi

\newsiamremark{remark}{Remark}
\newsiamremark{hypothesis}{Hypothesis}
\crefname{hypothesis}{Hypothesis}{Hypotheses}
\newsiamthm{claim}{Claim}

\headers{Tensor Completion via Gaussian Process based initialization}{Y. Kapushev, I.Oseledets, and E. Burnaev}

\title{Tensor Completion via Gaussian Process based initialization\thanks{Submitted to the editors DATE.
\funding{This  work  was  supported  by  the  Ministry  of  Education  and  Science  of  the  Russian  Federation  (Grant  no.14.756.31.0001)}}}

\author{Yermek Kapushev\thanks{Center for Computational and Data-Intensive Science and Engineering of Skolkovo Institute of Science and Technology, Moscow, Russia (\email{y.kapushev@skoltech.ru}, \email{i.oseledets@skoltech.ru}, \email{e.burnaev@skoltech.ru}).}
\and Ivan Oseledets\footnotemark[2]
\and Evgeny Burnaev\footnotemark[2]}

\usepackage{amsopn}

\usepackage{todonotes}
\usepackage{subcaption}

\ifpdf
\hypersetup{
  pdftitle={Tensor Completion via Gaussian Process based initialization},
  pdfauthor={Y. Kapushev, I. Oseledets, and E. Burnaev}
}
\fi

\begin{document}

\maketitle

\begin{abstract}
In this paper, we consider the tensor completion problem representing the solution in the tensor train (TT) format.
It is assumed that tensor is of high order, and tensor values are generated by an unknown smooth function.
The assumption allows us to develop an efficient initialization
scheme based on Gaussian Process Regression and TT-cross approximation technique.
The proposed approach can be used in conjunction with any optimization algorithm that is usually utilized in tensor completion problems.
We empirically justify that in this case the reconstruction error improves
compared to the tensor completion with random initialization.
As an additional benefit, our technique automatically selects rank thanks to using the TT-cross approximation technique.
\end{abstract}



\section{Introduction}

In this paper we consider the tensor completion problem.
We suppose that values of tensor $\mathcal{X}$ are generated by some smooth function, i.e.
\[
\mathcal{X}_{i_1, \ldots, i_d} = f(x_{i_1}, \ldots, x_{i_d}),
\]
where $(x_{i_1}, \ldots, x_{i_d})$ is a point on
some multi-dimensional grid and $f(\cdot)$ is some
unknown smooth function.
However, the tensor values are known only at some small subset of the grid.
The task is to complete the tensor, i.e., to reconstruct the tensor values at all points on the grid taking into account the properties of the {\em data generating process} $f(\cdot)$.

This problem statement differs from the traditional problem statement, which does not use any assumptions about the function $f(\cdot)$.
Knowing some properties of the data generating function provides insights about how the tensor values relate to each other, and this, in turn, allows us to improve the results.
In this work we assume that function $f(\cdot)$ is smooth.

There are a lot of practical applications that suit the statement.
For example, modeling of physical processes, solutions of differential equations, modeling probability distributions.

In this paper we propose to model the smoothness of the data generating process by using Gaussian Process Regression (GPR).
In GPR the assumptions about the function that we approximate are controlled via the kernel function.
The GPR model is then used to construct the initial solution to the tensor completion problem.

In principle, such initialization can improve any other tensor completion technique.
It means that using the proposed initialization state-of-the-art results
can be obtained employing some simple optimization procedure like Stochastic Gradient Descent.

When the tensor order is high the problem should be solved in some low-rank
format because the number of elements of the tensor grows exponentially.
The proposed approach is based on the tensor train (TT) format for its
computational efficiency and ability to handle large dimensions \cite{oseledets2010tt}.


The contributions of this paper are as follows
\begin{itemize}
\item[--] We introduce new initialization algorithm which takes into
account the tensor generating process.
The proposed algorithm is described in \cref{sec:initialization}.
\item[--] The proposed initialization technique automatically selects the rank of the tensor, the details are given in \cref{sec:tt_cross}.
\item[--] We conducted empirical evaluations of the proposed approach and compared it with tensor completion techniques without our initialization.
The results are given in \cref{sec:experiments} and show the superiority of the proposed algorithm.
\end{itemize}

\section{Tensor completion}
\label{sec:tensor_completion}

The formal problem statement is as follows.
Suppose that $\mathcal{Y}$ is a $d$-way tensor,
$\mathcal{Y} \in \mathbb{R}^{n_1 \times n_2 \times \cdots \times n_d}$
(by tensor here we mean a multi-dimensional array).
Tensor values are known only at some subset of indices
$\Omega \subset \{1, \ldots, n_1\} \times \cdots \times \{1, \ldots, n_d\}$.
By $P_\Omega$ we denote the projection onto the set $\Omega$, i.e.
\[
P_{\Omega} \mathcal{X} = \mathcal{Z}, \quad
\mathcal{Z}(i_1, i_2, \ldots, i_d) = \begin{cases}
\mathcal{X}(i_1, i_2, \ldots, i_d) & \mbox{if } (i_1, \ldots, i_d) \in \Omega, \\
0 & \mbox{otherwise}.
\end{cases}
\]
We formulate the tensor completion as an optimization problem
\begin{equation}
\begin{aligned}
\label{eq:tensor_completion}
    \min_{\mathcal{X}} \quad &  f(\mathcal{X}) = \|P_{\Omega} \mathcal{X} - P_{\Omega} \mathcal{Y}\|_F^2 \\
    \mbox{subject to} \quad & \mathcal{X} \in \mathcal{M}_r =
    \{\mathcal{X} \in \mathbb{R}^{n_1 \times \cdots n_d} \; | \; {\rm rank}_{TT}
    (\mathcal{X}) = \mathbf{r} \},
\end{aligned}
\end{equation}
where ${\rm rank}_{TT}(\mathcal{X})$ is a tensor train rank of $\mathcal{X}$ \cite{oseledets2011tensor},
which is a generalization of the matrix rank,
and $\|\cdot\|_F$ is the Frobenius norm.
A tensor $\mathcal{X}$ is said to be in tensor train format if its elements are represented as 
\[
    \mathcal{X}(i_1, \ldots, i_d) = \sum_{j_1, j_2, \ldots j_d}
    \mathcal{G}^{(1)}_{1, i_1, j_1}
    \mathcal{G}^{(2)}_{j_1, i_2, j_2} \cdots
    \mathcal{G}^{(d)}_{j_{d - 1}, i_d, 1},
\]
where $\mathcal{G}^{(i)}$ is a three-way tensor core with size
$r_{i - 1} \times n_i \times r_{i}$, $r_0 = r_{d} = 1$.
Vector $\mathbf{r}_{TT} = (r_0, \ldots, r_d)$ is called TT-rank.

Tensor train format assumes that the full tensor can be approximated by
a set of $3$-way core tensors, the total number of elements in
core tensors is $\mathcal{O}(dnr^2)$, where
$r = \max\limits_{i = 0, \ldots, d}\{r_i\}$,
$n = \max\limits_{i = 1, \ldots, d}\{n_i\}$,
which is much smaller than $n^d$.

In problem \eqref{eq:tensor_completion} we optimize the objective
function straightforwardly with respect to tensor cores
$\mathcal{G}^{(1)}, \ldots \mathcal{G}^{(d)}$ while having their sizes fixed.
Problem \eqref{eq:tensor_completion} is non-convex, so
optimization methods can converge to a local minimum.
To get an efficient solution we impose two requirements:
\begin{enumerate}
    \item Initial tensor $\mathcal{X}_0$ in tensor train format should be as close to the optimum as possible.
    \item Availability of an efficient optimization procedure that will be launched from the obtained initial tensor.
\end{enumerate}
These steps are independent, and one can apply any desired algorithm in each of them.

In this work we develop the initialization algorithm that allows obtaining accurate initial tensor for the case when the tensor of interest is generated by some smooth function.
The experimental section demonstrates that our initialization can improve the results of many optimization procedures,
and shows the potential of our approach to be adapted to a large number of different tensor completion techniques.


\section{Initialization}
\label{sec:initialization}
We consider tensors that are generated by some function,
i.e. tensor values are computed as follows
\[
\mathcal{Y}_{i_1, \ldots, i_d} = f(x_{i_1}, \ldots, x_{i_d}),
\]
where $f(\cdot)$ is some unknown smooth function and
$(x_{i_1}, \ldots, x_{i_d})^\top \in \mathbb{R}^d$,
$i_k = 1, \ldots, n_k$,
$n_1, \ldots, n_d$ are tensor sizes.
The set of points $\{(x_{i_1}, \ldots x_{i_d}): i_k = 1, \ldots, n_k; k = 1, \ldots, d\}$ is a full factorial Design of Experiments, i.e. a multi-dimensional grid,
and we also assume that the grid is uniform.

In this setting the tensor completion can be considered as a regression problem and can be solved by any regression technique that guarantees the smoothness of the solution.
However, in the tensor completion problem we are interested in a tensor of values of $f(\cdot)$ at a predefined finite grid of points.
The tensor should be in a low-rank format to be able to perform various operations with tensor efficiently (e.g. calculation of the norm of the tensor, dot product and other).

These observations give us the solution --- build regression model $\widehat{f}$ using the observed values of the tensor, then use the obtained approximation as a black-box for the TT-cross approximation algorithm \cite{oseledets2010tt}.
The last step results in a tensor $\widehat{\mathcal{X}}$ in TT format, which is a low-rank format and allows efficient computations.
The next step (which is optional) is to improve the obtained solution $\widehat{\mathcal{X}}$ by using it as initialization for any other tensor completion technique.


Let us write down the set of observed tensor values into a vector $\mathbf{y}$
and the corresponding indices into a matrix $\mathbf{X}$ (each row is a vector of indices
$(i_1, i_2, \ldots, i_d)$).
Then the approach for tensor completion (in TT format) can be
written as follows

\begin{enumerate}
    \item Construct initial tensor $\mathcal{X}_0$ in TT format:
        \begin{enumerate}
        \item Apply some regression technique using given data set $(\mathbf{X}, \mathbf{y})$ to construct approximation of the function that generates tensor values.
        \item Apply TT-cross method (see \cref{sec:tt_cross}, \cite{oseledets2010tt}) to the constructed approximation to obtain $\mathcal{X}_0$.
        \end{enumerate}
    \item Apply some tensor completion technique using $\mathcal{X}_0$
    as an initial value.
\end{enumerate}

At step $1$(a) the choice of the regression technique affects
the result of the initialization, although it can be arbitrary.
It is required to choose the regression algorithm such that it will capture the peculiarities of the tensor we would like to restore.
In this work we suppose that the tensor generating function is smooth
(which is a common situation when modeling physical processes).
Therefore, we choose a regression technique that is good at approximating smooth functions.
A reasonable choice, in this case, is to use Gaussian Process (GP) Regression 
\cite{rasmussen2004gaussian}.
GP models is a favorite tool in many engineering applications as they
proved to be efficient, especially for problems where it is required to
model some smooth function \cite{belyaev2016gtapprox}.
The points $(x_{i_1}, \ldots, x_{i_d})$ are not given, all we know is that at the point with multi-index
$(i_1, \ldots, i_d)$ on the grid the function value
is equal to $\mathcal{X}_{i_1, \ldots, i_d}$.
To make the problem statement reasonable we assume that
the indices are connected with the points as follows:
$x_{i_k} = a_k i_k + b_k$,
where $a_k, b_k \in \mathbb{R}$.
So, as an input for the approximation we set $a_k$ and $b_k$ such that $x_{i_k} \in [0, 1]$.

At step $1$(b) we use TT-cross
because it allows to efficiently approximate black-box function by a low-rank tensor in TT format.
Moreover, this approach can automatically select TT-rank making it more desirable.
More details on the technique are given in \cref{sec:tt_cross}.

The described approach has the following benefits:
\begin{enumerate}
    \item Initial tensor $\mathcal{X}_0$ which is close to the optimal value in terms of the reconstruction error at observed values.
    It will push the optimization to faster convergence.
    \item Better generalization ability --- there are many degrees of freedom: a lot of different tensor train factors can give low reconstruction error at observed positions but can give a large error at other locations.
    Accurate approximation model will push the initial tensor to be closer to the original tensor in both the observed positions and unobserved ones.
    \item TT-cross technique chooses rank automatically, so there is no need to tune the rank of the tensor manually.
\end{enumerate}

The described approach leads to the \cref{alg:initialization}.
The steps 3 and 4 of the algorithm are described in
\cref{sec:gp} and \cref{sec:tt_cross} correspondingly.

\begin{algorithm}
\caption{Initialization}
\label{alg:initialization}
\begin{algorithmic}
\STATE \textbf{Input}: $\mathbf{y}, \Omega$
\STATE{\textbf{Output}: $\mathcal{Y}_0$ in tensor train format}
\begin{enumerate}
    \item Construct the training set $(\mathbf{X}, \mathbf{y})$ from $\mathbf{y}, \Omega$
    \item Rescale inputs $\mathbf{X}$ to $[0, 1]$ interval
    \item
    \parbox[t]{0.5\textwidth}{
        Using $(\mathbf{X}, \mathbf{y})$ build GP model $\hat{f}(\mathbf{x})$
    }
    \hfill
    \parbox{0.8em}{$\triangleright$}
    \parbox[t]{0.3\textwidth}{
        see \cref{sec:gp} for details
    }
    \item \parbox[t]{0.5\textwidth}{
        Apply TT-cross to $\hat{f}(\mathbf{x})$ and obtain $\mathcal{Y}_0$
    }
    \hfill\parbox{0.8em}{$\triangleright$}
    \parbox[t]{0.3\textwidth}{
        see \cref{sec:tt_cross} for details
    }
\end{enumerate}
\RETURN{$\mathcal{Y}_0$}
\end{algorithmic}
\end{algorithm}

\subsection{Gaussian Process Regression}
\label{sec:gp}
One of the most efficient tools for approximating smooth functions is the Gaussian Process (GP)
Regression \cite{burnaev2016regression}.
GP regression is a Bayesian approach where a prior distribution over continuous functions
is assumed to be a Gaussian Process, i.e.
\[
  \mathbf{y} \, | \, \mathbf{X} \sim \mathcal{N}(\boldsymbol{\mu}, \, \mathbf{K}_f + \sigma_{noise}^2\mathbf{I}),
\]
where $\mathbf{y} = (y_1, y_2, \ldots, y_N)$ is a vector of outputs,
$\mathbf{X} = (\mathbf{x}_1^{\top}, \mathbf{x}_2^{\top}, \ldots, \mathbf{x}_N^{\top})^{\top}$ is a matrix of inputs,
$\mathbf{x}_i \in \mathbb{R}^d$,
$\sigma_{noise}^2$ is a noise variance,
$\boldsymbol{\mu} = (\mu(\mathbf{x}_1), \mu(\mathbf{x}_2), \ldots, \mu(\mathbf{x}_N))$ is a mean vector modeled by some function $\mu(\mathbf{x})$,
$\mathbf{K}_f = \{ k(\mathbf{x}_i, \mathbf{x}_j) \}_{i, j = 1}^N$ is a covariance matrix for some a priori selected covariance function $k$ and
$\mathbf{I}$ is an identity matrix.
An example of such function is a squared exponential kernel
\[
k(\mathbf{x}, \mathbf{x}') = \exp \left (
- \frac{1}{2} \sum_{i=1}^d \left (
\frac{\mathbf{x}^{(i)} - \mathbf{x}'{}^{(i)}}{\sigma_i}
\right )^2
\right ),
\]
where $\sigma_i, i = 1, \ldots, d$ are parameters of the kernel
(hyperparameters of the GP model).
The hyperparameters should be chosen according to the given data set.

Without loss of generality we make the standard assumption of zero-mean data.
Now, for a new unseen data point $\mathbf{x}_*$ we have
\begin{equation}
\label{eq:gp_posterior}
    \hat{f}(\mathbf{x}_*) \sim \mathcal{N}\left (\mu(\mathbf{x}_*), \sigma^2(\mathbf{x}_*) \right ),
\end{equation}
\[
    \mu(\mathbf{x}_*) = \mathbf{k}(\mathbf{x}_*)^\top \mathbf{K}_y^{-1}\mathbf{y},
\]
\[
    \sigma^2(\mathbf{x}_*) = k(\mathbf{x}_*, \mathbf{x}_*) - \mathbf{k}(\mathbf{x}_*)^\top \mathbf{K}_y^{-1} \mathbf{k}(\mathbf{x}_*),
\]
where $\mathbf{k}(\mathbf{x}_*) = (k(\mathbf{x}_*, \mathbf{x}_1), \ldots, k(\mathbf{x}_*, \mathbf{x}_N))^T$ and
$\mathbf{K}_y = \mathbf{K}_f + \sigma_{noise}^2\mathbf{I}$.


Let us denote the vector of hyperparameters $\sigma_i, i=1, \ldots, d, \sigma_f$ and $\sigma_{noise}$ by $\boldsymbol{\theta}$. 
To choose the hyperparameters of our model we consider the log-likelihood
\[
    \log p(\mathbf{y} \, | \, \mathbf{X}, \boldsymbol{\theta}) =
    -\frac12 \mathbf{y}^T \mathbf{K}_y^{-1}\mathbf{y} - \frac12 \log |\mathbf{K}_y| - 
    \frac{N}{2} \log 2 \pi
\]
and maximize it over the hyperparameters \cite{rasmussen2004gaussian}.
The runtime complexity of learning GP regression is $\mathcal{O}(N^3)$ as we need to calculate the inverse of $\mathbf{K}_y$, its determinant and derivatives of the log-likelihood.
If the sample size $N$ ($|\Omega|$ in our case) is large, the computational complexity becomes an
issue.
There are several ways to overcome it.
If the data set has a factorial structure (multidimensional grid in a simple
case) we can use the algorithm from \cite{belyaev2015gaussian}.
If the structure is factorial with a small number of missing values
the method from \cite{belyaev2016computationally} should be applied.
For general unstructured cases, the approximate GP model can be built
using, for example, the model described in \cite{munkhoeva2018quadrature} or use a subsample as a training set.

After tuning of the hyperparameters, we can use the mean of the posterior distribution
\eqref{eq:gp_posterior} as a prediction of the model.

Note, that the input points $\mathbf{X}$ in our case is a set of indices of the
observed values $\mathbf{y}$.
For the GP model we scale each index to $[0, 1]$ interval.

\subsubsection{Kernel choice and smoothness assumption}
We approximate function $f$ using the GP model.
The GP model is a function from some reproducing kernel Hilbert space (RKHS) $\mathcal{H}$ which is fully defined by the kernel function.
In ideal case the function $f$ should be from the Hilbert space $\mathcal{H}$.
However, for a given kernel function it can be difficult to identify what functions lie in the corresponding RKHS.

In practice GP models with popular kernels (RBF kernel, Mat{\'{e}}rn kernel) provide good results, when $f \in C^k$, $k \ge 1$.
So, in the experiments section we assume that the functions are from this class of functions and use RBF kernel.

\subsection{Tensor-Train cross-approximation}
\label{sec:tt_cross}
To approximate tensor $\widehat{\mathcal{X}}$ generated by $\hat{f}$ we use Tensor-Train cross-approximation.
First, let us consider the matrix case.
Suppose that we are given a rank-$r$ matrix $\mathbf{A}$ of size $m \times n$.
A cross-approxi\-mation for the matrix is represented as 
\[
    \mathbf{A} = \mathbf{C} \widehat{\mathbf{A}}^{-1} \mathbf{R},
\]
where $\mathbf{C} = \mathbf{A}(:, J), \mathbf{R} = \mathbf{A}(I, :)$ are some $r$ columns and rows of the matrix $\mathbf{A}$ and $\widehat{\mathbf{A}} = A(I, J)$ is the submatrix on the intersection of these rows and columns. To construct accurate approximation it is required to find submatrix $\widehat{\mathbf{A}}$ of large volume.
It can be done in $\mathcal{O}(nr^2)$ operations \cite{tyrtyshnikov2000incomplete}.

Now for tensor $\widehat{\mathcal{X}} \in \mathbb{R}^{n_1 \times \cdots \times n_d}$ the procedure is the following.
At the first step let us consider unfolding $\mathbf{X}_1$ of size
$n_1 \times n_2 n_3 \cdot \cdots \cdot n_d$ and rank $r_1$.
Using row-column alternating algorithm from \cite{tyrtyshnikov2000incomplete} we can find $r_1$ linearly independent columns of matrix $\mathbf{X}_1$, these columns form matrix $\mathbf{C}$.
After that applying maxvol procedure \cite{tyrtyshnikov2000incomplete} to the matrix $\mathbf{C}$ we can find set of row indices $I_1 = \big[i_1^{\alpha_1} \big], \alpha_1 = 1, \ldots, r_1$, matrix $\mathbf{R}$ and matrix $\widehat{\mathbf{A}}_1$
that will give the cross-approximation of unfolding $\mathbf{X}_1$:
\[
    \mathbf{X}_1 = \mathbf{C}\widehat{\mathbf{A}}_1^{-1}\mathbf{R}.
\]
We set
\[
    \mathbf{G}_1 = \mathbf{C}\widehat{\mathbf{A}}_1^{-1},
\]
where $\mathbf{G}_1$ is of size $n_1 \times r_1$.
Next, let us form tensor $\mathcal{R}$ from $r_1$ rows of $\mathbf{X}_1$:
\[
    \mathcal{R}(\alpha_1, i_2, \ldots, i_d) = \widehat{\mathcal{X}}(i_1^{\alpha_1}, i_2, \ldots, i_d),
\]
and reshape it into a tensor of size $r_1n_2 \times n_3 \times \cdots \times n_d$.
Next step is to apply the same procedure to the unfolding $\mathbf{R}_1$ of the tensor $\mathcal{R}$ and obtain the matrices $\mathbf{C}$, $\widehat{\mathbf{A}}_2$ and
\[
    \mathbf{G}_2 = \mathbf{C} \widehat{\mathbf{A}}_2^{-1}
\]
of size $r_1n_2 \times r_2$.

Repeating the described procedure $d$ times we will end up with matrices
$\mathbf{G}_1, \allowbreak \mathbf{G}_2, \ldots, \mathbf{G}_d$ of sizes $n_1 \times r_1, r_1n_2 \times r_2, \ldots, r_{d-1}n_d \times 1$.
Then each matrix can be reshaped to the $3$-way tensor of size $r_{d - 1} \times n_d \times r_d$,
$r_0 = r_d = 1$ and can be used as core tensors for TT format.
It turns out that such representation is a TT decomposition of the initial tensor $\widehat{\mathcal{X}}$.

The exact ranks $r_1, \ldots, r_d$ are not known to us in general.
They can only be estimated from the above (e.g., by the maximum rank of the corresponding unfolding).
If the rank is overestimated then the calculation of matrices $\mathbf{G}_i$ 
is an unstable operation (because we obtain almost rank-deficient unfolding matrices).
However, in \cite{oseledets2010tt} the authors suggest some simple modification that overcomes this issue.
Therefore, we need to estimate the ranks from the above, but
the estimate should not be much larger than the real rank.
So, the approach is to start from some small rank, construct the tensor
in TT format and then apply recompression (see \cite{oseledets2011tensor}).
If there is a rank that is not reduced, then we underestimated that rank and should increase it and repeat the procedure.

\subsection{Computational complexity}
The computational complexity of TT cross-approximation method is as follows.
To perform the procedure we need to evaluate the GP model $\mathcal{O}(dnr^2)$ times at some subset of grid points
and then perform $\mathcal{O}(dnr^3)$ operations to find all maximum volume submatrices.
The complexity of evaluating GP model at one point is $\mathcal{O}(Nd)$, where $N$ is the number of observed tensor elements.
The total complexity is thus $\mathcal{O}(Nd^2nr^2 + dnr^3)$.

\section{Experimental results}
\label{sec:experiments}
In this section we present the results of the application of our approach
to two engineering problems and also test it on some artificial problems to investigate how its properties depend on smoothness.

The experimental setup is the following.
We try the following optimization algorithms
\begin{enumerate}
    \item SGD -- stochastic gradient descent \cite{zhao2018high},
    \item Ropt -- Riemannian optimization \cite{steinlechner2016riemannian},
    \item TTWopt -- weighted tensor train optimization \cite{zhao2018high},
    \item ALS  -- alternating least squares \cite{grasedyck2013alternating}.
\end{enumerate}
We run each algorithm with random initialization and with the
proposed GP-based initialization and then compare results.

\subsection{Functions generated from GP prior}
In order to study the dependence of the solution on the smoothness of the generating function we applied the proposed approach to the toy functions generated from GP prior with different kernels.
The smoothness of the generated functions is the same as the kernel that we used to generate them.
Thus, we can investigate the performance of the approach for different smoothness.
We considered shift-invariant kernels, i.e. $k(\mathbf{x}, \mathbf{y}) = k(r)$, where $r = \|\mathbf{x} - \mathbf{y}\|$.
The list of kernels is as follows:
\begin{itemize}
    \item Exponential kernel
    \[
        k(r) = \exp \left ( -\frac{r}{\sigma} \right ).
    \]
    The kernel is not differentiable at $\mathbf{x} = \mathbf{y} = 0$, thus
    the functions generated with this kernel are from $C^0$.
    \item Matern$_{3/2}$
    \[
        k_{3/2}(r) = \left (1 + \frac{\sqrt{3}r}{\sigma} \right )
        \exp \left ( -\frac{\sqrt{3}r}{\sigma} \right).
    \]
    This kernel is $1$-time differentiable.
    \item Matern$_{5/2}$
    \[
        k_{5/2}(r) = \left (1 + \frac{\sqrt{5}r}{\sigma} + \frac{5r^2}{3\sigma^2}\right )
        \exp \left ( -\frac{\sqrt{5}r}{\sigma} \right).
    \]
    This kernel is $2$-times differentiable.
    \item Radial Basis Function (RBF) kernel
    \[
        k(r) = \exp \left ( -\frac{r^2}{2\sigma^2} \right ).
    \]
    This kernel is infinitely differentiable.
    
\end{itemize}
Note, that despite the functions were generated using different kernel functions
in the proposed approach we used RBF kernel.

To compare how much one approach is better than the other we calculate the relative MSE error
\[
    {\rm MSE_{rel}} = \frac{1}{\lvert \Omega_{test}\rvert} \left \|\frac{P_{\Omega_{test}}\hat{\mathcal{Y}} -
      P_{\Omega_{test}}\mathcal{Y}}
      {\hat{\sigma}}\right \|_F^2,
\]
where $\Omega_{test}$ is some set of indices independent from the given observed set of indices $\Omega$, $|\Omega_{test}|$ is a size of the set $\Omega_{test}$ and $\widehat{\mathcal{Y}}$ is an obtained approximation of the actual tensor $\mathcal{Y}$,
$\hat{\sigma}$ is a standard deviation of $P_{\Omega_{test}}\mathcal{Y}$.
Such error can be interpreted as the ratio of unexplained variance.
For each optimization technique we calculate the difference between error obtained using random initialization and the error obtained using GP based initialization.

\begin{figure}
    \centering
    \includegraphics[width=\textwidth]{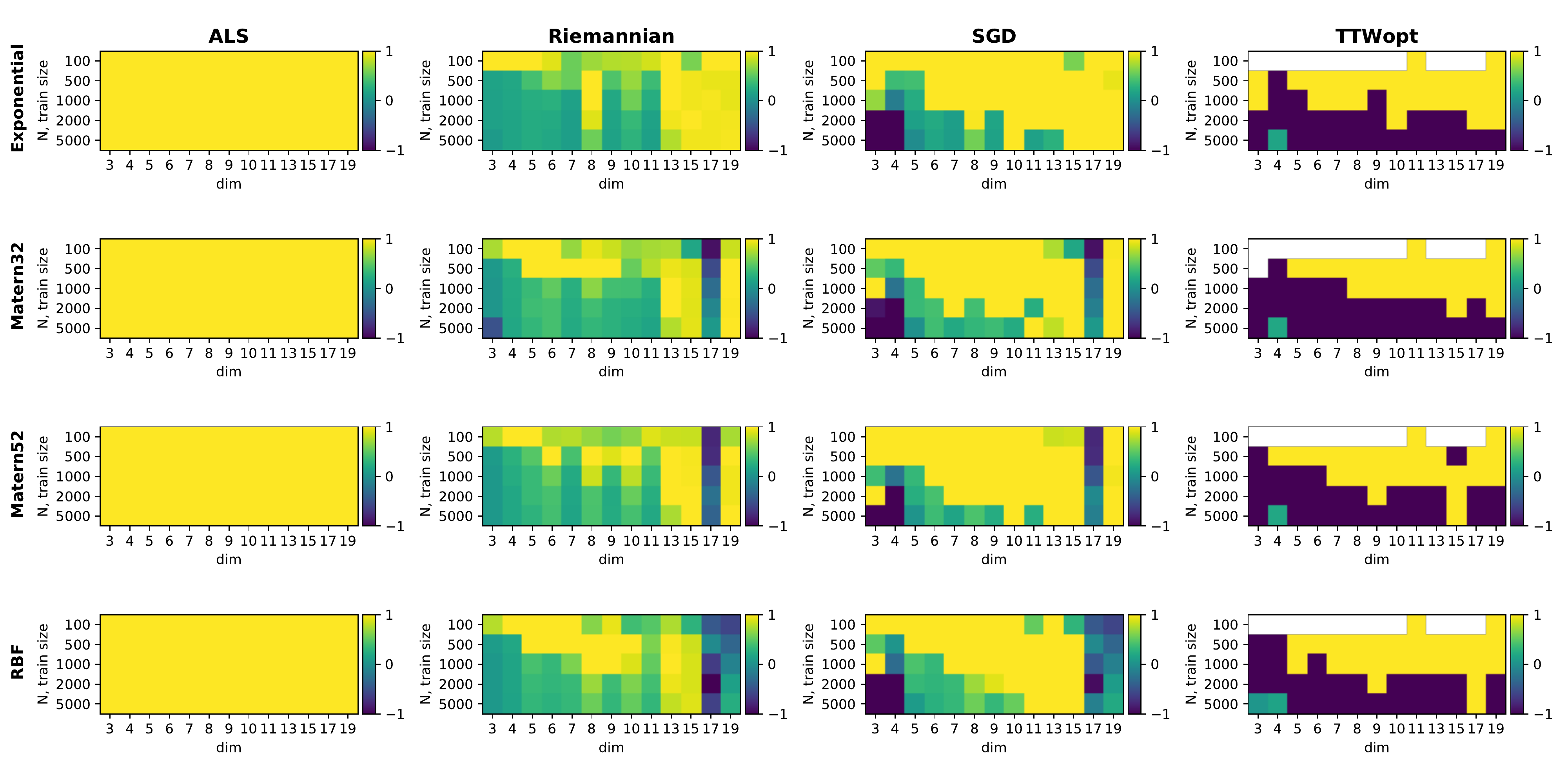}
    \caption{Improvement of GP based initialization over random initialization for different tensors and optimization methods. We clamp the improvement value to $[-1, 1]$ interval to make plots more illustrative.}
    \label{fig:improvement}
\end{figure}

For each kernel function we generated several data sets with different number of observed points ($N \in \{100, 500, 1000, 2000, 5000\}$) and different dimensionalities ($d \in \{2, 3, \ldots, 10, 11, 13, \ldots, 19\} $).
The quantity is illustrated in Figure \ref{fig:improvement} (note, that we clamp the values to $[-1, 1]$ interval to make the figures more illustrative).
It shows the improvement of one initialization over another.
We can see from the figure that in most of the cases GP-based initialization gives high improvement in the relative error.
The only exception is TTWopt method for which the benefit of the proposed initialization scheme takes place only in about half of the cases.
You can also note white squares for TTWopt, they mean that implementation of TTWopt crashed for the given data set (for some unknown reason).

\subsection{Read world functions}
We compared the approaches on two real world problems:
CMOS oscillator model and Cookie problem (see \cref{sec:cmos_ring_oscillator} and \cref{sec:cookie_problem} correspondingly).
In CMOS oscillator problem we run each optimization $10$ times with different random training sets and then calculate the average reconstruction error as well as standard deviation.
Cookie problem is more computationally intensive because each evaluation of the tensor value takes more resources and time.
Therefore, for Cookie problem we performed $10$ runs, and the training set was the same during all runs.

The quality of the methods is measured using mean squared error (MSE)
\[
    MSE = \frac{1}{|\Omega_{test}|}\|P_{\Omega_{test}} \widehat{\mathcal{Y}} - P_{\Omega_{test}}\mathcal{Y} \|_F^2.
\]
We also report the error of the initial tensor for random and the proposed initializations for each problem.

Note, that when we use GP based initialization, the TT-rank
$\mathbf{r}_{{\rm TT}}$ of the tensor is selected automatically by the TT-cross algorithm and max value of $\mathbf{r}_{{\rm TT}}$ can be larger than $n$.
The optimization algorithms with random initialization do not have a procedure for automatic rank selection, so we ran them with different ranks (from $1$ to $\min_{k}{n_k}$) and then chose the best one.

TTWopt implementation\footnote{\url{https://github.com/yuanlonghao/T3C_tensor_completion}} does not support high-dimensional problems.
For higher dimensional problems the authors of TTWopt propose to use SGD.
The authors of TTWopt also propose truncated SVD based initialization.
The idea is to fill missing values using the mean value of the
observed part of the tensor and then to apply truncated SVD
to obtain TT cores.
However, such approach is only applicable to low-dimensional
tensors as it requires to calculate full matrices of large
size.

For Ropt and ALS we used publically available MATLAB codes
\footnote{\url{https://anchp.epfl.ch/index-html/software/ttemps/}}.

\subsection{Cookie problem}
\label{sec:cookie_problem}
\begin{figure}[htbp]
  \centering
  \includegraphics[width=0.5\textwidth]{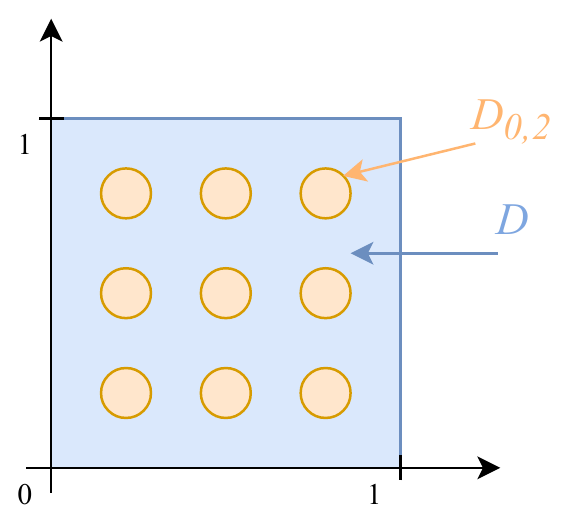}
  \caption{Illustration of Cookie problem with $m=3$ ($9$ cookies).}
  \label{fig:cookie_problem}
\end{figure}

Let us consider parameter-dependent PDE
\cite{ballani2015hierarchical, tobler2012low}:
\begin{align*}
    -{\rm div} (a(x, p) \nabla u(x, p)) &= 1, \quad x \in D = [0, 1]^2, \\
    u(x, p) &= 0, \quad x \in \partial D,
\end{align*}
where
\[
a(x, p) = \begin{cases}
p_\mu, \quad \mbox{if } x\in D_{s, t}, \mu = mt + s, \\
1, \quad \mbox{otherwise},
\end{cases}
\]
$D_{s, t}$ is a disk of radius $\rho=\frac{1}{4m + 2}$
and $m^2$ is a number of disks which form $m \times m$ grid.
This is a heat equation where heat conductivity $a(x, p)$ depends on $x$
(see illustration in \cref{fig:cookie_problem}) and $p$ is an $m^2$-dimensional.

We are interested in average temperature over $D$:
$u(p) = \int_{[0, 1]^2} u(x, p){\rm d}x$.
If $p$ takes $10$ possible values then there are $\mathbf{10^{m^2}}$
possible values of $u(p)$.

In this work we used the following setup for the Cookie problem: each parameter $p$ lie in the interval
$[0.01, 1]$,
number of levels for each p is $10$,
number of cookies is $m^2 = 9$ and $16$,
size of the observed set is $N = 5000$,
for the test set we used $10000$ independently generated points.

The results of tensor completion are presented in \cref{tab:results_cookie} (the variance of the initialization error for GP-based init is not presented as it is negligible in this case). One can see that GP based initialization gives lower reconstruction errors both on the training set and test set except for ALS technique. ALS method with the proposed initialization
overfits: the error on the training set is close to $0$, whereas the test error is much more significant. The error on the training set is about $10^{-29}$, which means that the training set was approximated with machine precision.
It is not surprising if we recall that there are only $5000$ observed values, while the number of free parameters that are used to construct TT is much higher.

\subsection{CMOS ring oscillator}
\label{sec:cmos_ring_oscillator}
Let us consider the CMOS ring oscillator \cite{zhang2017big}. It is an electric circuit which consists of $7$ stages of CMOS inverters.
We are interested in the oscillation frequency of the oscillator. The characteristics of the electric circuit are described by $57$ parameters.
Each parameter can take one of $3$ values, so the total size of the tensor is $3^{57} \approx 1.57 \times 10^{27}$. The number of observed values that were used during the experiments is $N = 5000$. For the test set we used $10000$ independently generated points.

The results of the experiments are given in \cref{tab:results_cmos}. The table demonstrates that utilizing GP based initialization improves the results for all algorithms except ALS. ALS, in this case, overfits again: training error is extremely small, whereas the test error is much larger, though it is rather small compared to other techniques and ALS with random initialization.

\begin{table}[]
    \centering
    \caption{MSE errors for Cookie problem}
    \small
    \label{tab:results_cookie}
    \begin{tabular}{r|cc|cc|}
         & \multicolumn{4}{c|}{\textbf{Training set}} \\
         & \multicolumn{4}{c|}{$m = 3$} \\
        \hline
          & \multicolumn{2}{c|}{Random init} & \multicolumn{2}{c|}{GP init} \\
          & \multirow{2}{*}{error} & \multirow{2}{4em}{average N iters} & \multirow{2}{*}{error} &
          \multirow{2}{4em}{average N iters} \\
          & & & & \\
        \hline
        SGD & $(1.66 \pm 0.067) \times 10^{-2}$ & $1500$ &
        $\mathbf{(2.86 \pm 0.18) \times 10^{-5}}$ & $150$ \\
        Ropt & $(4.13 \pm 2.20) \times 10^{-8}$ & $1000$ &
        $\mathbf{(5.48 \pm 1.10) \times 10^{-10}}$ & $1000$ \\
        TTWopt & $(2.73 \pm 0.19) \times 10^{-4}$ & $100$ &
        $\mathbf{(9.21 \pm 2.17) \times 10^{-7}}$ & $100$ \\
        ALS & $(1.07 \pm 1.07) \times 10^{-4}$ & $100$ &
        $\mathbf{(2.39 \pm 0.60) \times 10^{-30}}$ & $100$ \\
        \hline
        Init error & $(1.15 \pm 0.24) \times 10^{8}$ &  &
        $\mathbf{2.66 \times 10^{-4}}$ &  \\
        \hline
        
        & \multicolumn{4}{c|}{$m = 4$} \\
        
        \hline
        SGD & 
        $(3.14 \pm 1.08) \times 10^{-2}$ & $1500$ & 
        $\mathbf{(1.65 \pm 0.13) \times 10^{-5}}$ & $150$ \\
        Ropt &
        $(1.42 \pm 0.01) \times 10^{-2}$ & $1000$ & 
        $\mathbf{(3.42 \pm 0.50) \times 10^{-4}}$ & $1000$ \\
        TTWopt &
        $(1.31 \pm 0.00) \times 10^{-4}$ & $100$ & 
        $\mathbf{(1.80 \pm 0.16) \times 10^{-6}}$ & $100$ \\
        ALS &
        $(6.59 \pm 3.30) \times 10^{-5}$ & $100$ & 
        $\mathbf{(1.33 \pm 0.46) \times 10^{-29}}$ & $100$ \\
        \hline
        Init error & $(8.32 \pm 2.52) \times 10^{14}$ &  &
        $\mathbf{3.14 \times 10^{-4}}$ &  \\
        \hline
        
         & \multicolumn{4}{c|}{\textbf{Test set}} \\
         & \multicolumn{4}{c|}{$m = 3$} \\
         \hline 
         
        SGD & $(2.06 \pm 2.31) \times 10^{-1}$ & --- & 
        $\mathbf{(9.97 \pm 0.40) \times 10^{-5}}$ & --- \\
        Ropt & $\mathbf{(1.48 \pm 0.90) \times 10^{-7}}$ & --- & 
        $(3.45 \pm 0.0165) \times 10^{-4}$ & --- \\
        TTWopt & $(4.52 \pm 0.50) \times 10^{-4}$ & --- & 
        $\mathbf{(5.27 \pm 0.74) \times 10^{-6}}$& --- \\
        ALS & $\mathbf{(4.37 \pm 7.73) \times 10^{-2})}$ & --- & 
        $(3.78 \pm 1.08) \times 10^{0}$ & --- \\
        \hline
        Init error & $(1.12 \pm 0.23) \times 10^{8}$ &  &
        $\mathbf{4.12 \times 10^{-4}}$ &  \\
        \hline

        & \multicolumn{4}{c|}{$m = 4$} \\
        
        \hline
        SGD & $(2.40 \pm 2.76) \times 10^{1}$ & --- & 
        $\mathbf{(1.15 \pm 0.05) \times 10^{-4}}$ & ---\\
        Ropt & $(1.47 \pm 0.003 \times 10^{-2}$ & --- &
        $\mathbf{(5.38 \pm 0.07) \times 10^{-4}}$ & ---\\
        TTWopt & $(2.42 \pm 0.00) \times 10^{-4}$ & --- &
        $\mathbf{(3.02 \pm 0.17) \times 10^{-5}}$ & --- \\
        ALS & $\mathbf{(3.57 \pm 5.65) \times 10^{-1}}$ & --- &
        $(1.85 \pm 60.5) \times 10^{0}$ & ---\\
        \hline
        Init error & $(8.33 \pm 2.46) \times 10^{14}$ &  &
        $\mathbf{5.37 \times 10^{-4}}$ &  \\
        \hline

    \end{tabular}
\end{table}

\begin{table}[]
    \centering
    \caption{MSE errors for CMOS oscillator}
    \label{tab:results_cmos}
    \begin{tabular}{r|cc|cc|}
         & \multicolumn{4}{c|}{\textbf{Training set}} \\
        \hline
        & \multicolumn{2}{c|}{Random init} & \multicolumn{2}{c|}{GP init} \\
        & error & N iters & error & N iters \\
        \hline
        SGD & $(7.77 \pm 15.25) \times 10^5$ & $1500$ &
        $\mathbf{(3.11 \pm 4.87) \times 10^{-4}}$ & $150$ \\

        Ropt & $(6.22 \pm 0.01) \times 10^3$ & $1000$ &
        $\mathbf{(9.50 \pm 4.28) \times 10^{-5}}$ & $1000$ \\
         
        ALS & $(9.95 \pm 0.26) \times 10^{-2}$ & $300$ & 
        $\mathbf{(3.57 \pm 0.45) \times 10^{-26}}$ & $300$ \\
        \hline
        Init error & $(1.67 \pm 3.19) \times 10^{12}$ & & 
        $\mathbf{(3.95 \pm 2.19) \times 10^{-4}}$ & \\
        \hline
        
         & \multicolumn{4}{c|}{\textbf{Test set}} \\
        \hline
        
        SGD & $(3.45 \pm 9.68) \times 10^8$ & --- &
        $\mathbf{(4.65 \pm 5.01) \times 10^{-4}}$ & --- \\
         
        Ropt & $(6.23 \pm 0.0) \times 10^3$ & --- & 
        $\mathbf{(9.68 \pm 4.16) \times 10^{-5}}$ & --- \\
         
        ALS & $(1.03 \pm 0.01) \times 10^{-1}$ & --- &
        $\mathbf{(4.09 \pm 3.10) \times 10^{-4}}$ & --- \\
        \hline
        Init error & $(1.04 \pm 2.53) \times 10^{15}$ & & 
        $\mathbf{(3.90 \pm 2.15) \times 10^{-4}}$ & \\
        \hline        
    \end{tabular}
\end{table}

All in all, the obtained results prove that GP based initialization allows improving the tensor completion results in general. At least it provides better training error. As for the error on the test set one should be more careful as the number of degrees of freedom is large and there exist many solutions that give a small error for the observed values but large errors for other values.

\section{Related works}
\label{sec:related_works}
One set of approaches to tensor completion is based on nuclear norm minimization. The nuclear norm of a matrix is defined as a sum of all singular values of the matrix. This objective function is a convex envelope of the rank function. For a tensor the nuclear norm is defined as a sum of singular values of matricizations of the tensor.

There are efficient off-the-shelf techniques for such types of problems that apply interior-point methods. However, they are second-order methods and scale poorly with the dimensionality of the problem. Special optimization technique was derived for nuclear norm minimization
\cite{gandy2011tensor, liu2013tensor, recht2010guaranteed}.



More often such techniques are applied to matrices or
low-dimensional tensors as their straightforward
formulation allows finding the full tensor.
It becomes infeasible when we come to high-dimensional problems.


The second type of approaches is based on low-rank tensor decomposition \cite{acar2011scalable, chen2013simultaneous, kressner2014low, steinlechner2016riemannian, yuan2017completion}. There are several tensor decompositions, and all these papers derive some optimization procedure for one of them, namely, CP decomposition, Tucker decomposition or TT/MPS decomposition. The simplest technique is the alternating least squares \cite{grasedyck2015alternating}. It just finds the solution iteratively at each iteration minimizing the objective function w.r.t. one core while other cores are fixed.

Another approach is based on Riemannian optimization that tries to find the optimal solution on the manifold of low-rank tensors of the given structure \cite{steinlechner2016riemannian}. The same can be done by using Stochastic Gradient Descent \cite{yuan2017completion}. Riemannian optimization, TTWopt, ALS and its modifications (e.g., ADF, alternating directions fitting \cite{grasedyck2013alternating}) try to find the TT representation of the actual tensor iteratively.  At each iteration it optimizes TT cores such that the resulting tensor approximates well the tensor which coincides with the real tensor at observed indices and with the result of the previous iteration at other indices. All these approaches need to specify rank manually. The authors of \cite{zhao2015bayesian} apply the Bayesian framework for CP decomposition which allows them to select the rank of the decomposition automatically.

In some papers the objective is modified by introducing special regularizers to suit the problem better \cite{yokota2016smooth}. For example, in \cite{chen2013simultaneous,zhao2015bayesian} to obtain better results for visual data a special prior regularizer was utilized.

Our proposed algorithm is an initialization technique for the
tensor completion problems in TT format and can be used
with most of the algorithms solving such problems.
If the assumptions from \cref{sec:initialization} (the tensor values are values of some rather smooth function of tensor indices) are satisfied the initial value will be close to the optimal
providing better results.
The question of a good initialization is rarely taken into account.
In paper \cite{ko2018fast} a special initialization if proposed for visual data. The idea is to use some crude technique (like bilinear interpolation) to fill missing values and after that apply SVD-based tensor train decomposition.
The drawback of the approach is that it can be applied only in case of small-dimensional tensors as we need to fill all missing values.
In \cite{grasedyck2013alternating} they propose special initialization for the Alternating Direction Fitting (ADF) method.
This is a general technique for the tensor completion and it does not take into account the assumptions on the data generating function.


\section{Conclusions}
\label{sec:conclusions}

We proposed a new initialization algorithm for high-dimen\-sional tensor completion in TT format. The approach is designed mostly for the cases when some smooth function generates the tensor values. It can be combined with any optimization procedure that is used for tensor completion. Additionally, the TT-rank of the initial tensor is adjusted automatically by the TT-cross method and defines the resulting rank of the tensor. So, the approach provides an automatic rank selection. Our experimental study confirms that the proposed initialization delivers lower reconstruction errors for many of the optimization procedures.

\section*{Acknowledgments}
We would like to thank Zheng Zhang, who kindly provided us CMOS ring oscillator data set.

\bibliographystyle{siamplain}
\bibliography{references}
\end{document}